\documentclass[11pt]{article}

\usepackage[centertags]{amsmath}
\usepackage{amsfonts}
\usepackage{amssymb}
\usepackage{theorem}
\usepackage{ulem}
\usepackage{epsfig}
\usepackage{subfigure}
\usepackage[all]{xy}
\usepackage{changebar}

\theoremstyle{break} \newtheorem{theorem}{Theorem}[section]
\theoremstyle{break} 
\theoremstyle{break}        
\theoremstyle{break} \newtheorem{lemma}[theorem]{Lemma}
\theoremstyle{break} \newtheorem{corollary}[theorem]{Corollary}
\theoremstyle{break} 
\theoremstyle{break} 
\theoremstyle{break}
{\theorembodyfont{\rmfamily}\newtheorem{remark}[theorem]{Remark}}
{\theorembodyfont{\rmfamily}}
\theoremstyle{break} 
\theoremstyle{break} 
\theoremstyle{break} 
\theoremstyle{break} 
\numberwithin{equation}{section}

\newcommand{\R}{{\mathbb{R}}}
\newcommand{\D}{{\mathbb{D}}}
\newcommand{\C}{{\mathbb{C}}}

\def\Re{\mathop{{\rm Re}}}

\begin{document}

\renewcommand{\thefootnote}{}
\stepcounter{footnote}
\begin{center}
{\bf \Large 
On the isolated singularities of the solutions 
of the\\[2mm] Gaussian curvature equation for nonnegative curvature
\footnote{
The first author was supported by a HWP scholarship; the second author
received partial support  from the German--Israeli Foundation (grant G--809--234.6/2003).}}
\end{center}

\medskip

\renewcommand{\thefootnote}{\arabic{footnote}}
\setcounter{footnote}{0}
\begin{center}
{\large  \bf Daniela Kraus and Oliver Roth}\\[2mm]
{\small Universit\"at W\"urzburg,
 Mathematisches Institut, \\[0.5mm]     
 D--97074 W\"urzburg, Germany\\[1mm]
dakraus@mathematik.uni-wuerzburg.de\\ roth@mathematik.uni-wuerzburg.de}\\[1mm]
\end{center}

\centerline{\today}

\renewcommand{\thefootnote}{\arabic{footnote}}

\medskip
\begin{center}
\begin{minipage}{13cm} {\bf Abstract.}   
The precise asymptotic behaviour of the solutions to the twodimensional curvature
equation $\Delta u=k(z)\,e^{2 u}$
with $e^{2 u} \in L^1$ for bounded nonnegative curvature functions $-k(z)$
near isolated singularities is obtained.
\end{minipage}
\end{center}

\section{Results}

The aim of this note is to classify the isolated singularities of the
real--valued solutions
with finite energy of the twodimensional curvature equation $\Delta u=k(z) \, e^{2
  u}$ in the case of {\it variable} nonnegative curvature.
We first consider  solutions $u \in C^2(\D \backslash \{ 0 \})$ of 
this equation, where
$\D:=\{z \in \C \, : \, |z|<1\}$ denotes the open unit disk in the complex plane.

\begin{theorem} \label{thm:main}
Let $k \in C(\D \backslash \{ 0 \})$  be a bounded and nonpositive function and
 $u \in C^2(\D \backslash \{ 0 \})$  a real--valued solution to the curvature
  equation
\begin{equation} \label{eq:gauss}
\Delta u= k(z) \, e^{2 u} 
\end{equation}
with
\begin{equation} \label{eq:int}
\iint \limits_{\D} e^{2 u(z)} \, dxdy <\infty \,  . 
\end{equation}
Then there exists a real number $\gamma >-1$ such that
\begin{equation} \label{eq:expansion}
u(z)=\gamma \log |z|+r(z) 
\end{equation}
where $r \in C(\D)$. In addition, 
 $r \in C^1(\D)$ if $\gamma>-1/2$, and 
$r \in
C^2(\D)$ if $\gamma \ge 0$ and if $k$ is locally H\"older continuous on $\D$.
\end{theorem}

\begin{remark}
\begin{itemize}
\item[(a)] For $k(z)\equiv -4$  the PDE (\ref{eq:gauss}) reduces to the 
{\it Liouville equation} $\Delta u=-4 \, e^{ 2u}$. 
In this case Theorem \ref{thm:main} was proved by Chou and Wan
 \cite{CW94,CW95} using complex analysis and Liouville's classical
representation formula \cite{Lio1853} for the solutions to $\Delta u=-4 \,
e^{2 u}$. In the variable curvature case some (nonsharp)
estimates for the solutions $u$ of (\ref{eq:gauss}) satisfying the energy
estimate (\ref{eq:int}) have recently been obtained by Yunyan \cite{Yun03}
using blow--up analysis and the moving plane method.
\item[(b)] If $k(z)$ is {\it strictly} positive (and bounded), then 
condition (\ref{eq:int}) is redundant. The asymptotic behaviour of the
solutions in this case was found by McOwen \cite{McO93} (see also Heins
\cite{Hei62}); the regularity
properties of the remainder function are studied in \cite{KR}.
We note that the isolated singularites of the solutions in the classcial
constant case 
$k(z)=+4$ were first described by J.~Nitsche \cite{Nit57} and later by Chou
and Wan \cite{CW94,CW95}.
\item[(c)] The ``energy condition'' (\ref{eq:int}) cannot be dispensed with.
This can already be seen in the constant case $k(z) = -4$. Here  one obtains
very badly behaved solutions $$u(z) =\log \left(\frac{|g'(z)|}{1+|g(z)|^2}\right)$$ at $z=0$ by choosing
$g$  meromorphic on $\D \backslash \{ 0 \}$ with an essential singularity at
$z=0$ (see \cite{CW94}).
Condition (\ref{eq:int}) seems to be a ``canonical'' assumption when dealing
with nonnegative curvature (see e.g.~Br\'ezis \& Merle \cite{BM}, Hang \&
Wang \cite{HW}, Jost, Wang \& Zhou \cite{JWZ}, Yunyan \cite{Yun03} and the
references cited therein).
\item[(d)] 
The punctured unit disk plays no special role
in Theorem \ref{thm:main}-- we can replace it by
any domain in the complex plane with an isolated boundary point and obtain a corresponding version of Theorem \ref{thm:main} in this more general situation.
This will become apparent from the proof of Theorem \ref{thm:main}.
 \end{itemize}
\end{remark}

Geometrically, every function $u$ of Theorem \ref{thm:main} gives rise to a
conformal Riemannian metric $\lambda(z) \, |dz|:=e^{u(z)} \, |dz|$ with a {\it conical
singularity of order} $\gamma>-1$ at $z=0$, i.e., $\lambda(z) \, |dz|=|z|^{\gamma} \, e^{r(z)}
\, |dz|$ and curvature $-k(z)$. Theorem 1.1 contains precise information how
the connection (sometimes also called
Pre--Schwarzian or Christoffel symbol) 
$$ \Gamma_{\lambda}(z):=2 \frac{\partial \log \lambda(z)}{\partial z} $$
of $\lambda(z) \, |dz|$
and its projective connection (see \cite{JWZ}) or
Schwarzian (see \cite{Min97})
$$ S_{\lambda}(z):=\frac{\partial \Gamma_{\lambda}(z)}{\partial z}-\frac{1}{2}
\, 
\Gamma_{\lambda}(z)^2= 2 \left[  \frac{\partial^2 \log \lambda(z)}{\partial
    z^2}- \left(  \frac{\partial \log \lambda(z)}{\partial z} \right)^2
\right] \, $$
behave at the conical singularity:

\begin{corollary} \label{cor}
Let $\lambda(z) \, |dz|$ be a conformal metric on $\D \backslash \{ 0\}$
with curvature $-k(z)$ for some nonpositive bounded continuous function
$k : \D\backslash \{0\} \to \R$ and $\iint_{\D} \lambda(z)^2 \, dxdy<\infty$. Then
\begin{itemize}
\item[(a)] $\lim \limits_{z \to 0} z \, \Gamma_{\lambda}(z)=\gamma$; and
\item[(b)] $\lim \limits_{z \to 0} z^2 \, S_{\lambda}(z)=-\gamma (2+\gamma)/2$
  if $k$ is locally H\"older continuous on $\D$, 
\end{itemize}
where $\gamma>-1$ is the order of the conical singularity of $\lambda(z) \,
|dz|$ at $z=0$.
\end{corollary}


\begin{remark} \label{rem:2}
The Schwarzian $S_{\lambda}$ of a conformal metric plays an important role
in particular for metrics with {\it constant} curvature. This classical constant
curvature case is intimately related to complex analysis. In fact, if
$\lambda (z) \, |dz|$ has constant curvature, then
 $S_{\lambda}$ is  a {\it holomorphic} function with isolated
singularities exactly at the isolated singularities of the metric $\lambda(z) \, |dz|$.
Corollary \ref{cor} shows that for constantly curved singular metrics with 
$\iint \lambda(z)^2 dxdy<\infty$ the Schwarzian $S_{\lambda}$ has 
a pole of order $2$ at the conical singularities. We refer to 
\cite{Bie12,Bie16,CW94,E,HW,Hei62,JWZ,Lio1853,MT2002,McO93,Min97,Nit57,Pic1893,Pic1905,Poi1898,Tro1992,Yam88}
for more information about conformal metrics with constant curvature and
conical singularities.
\end{remark}

\section{Proofs}

Theorem \ref{thm:main} follows from the following lemma, the Br\'ezis--Merle
lemma \cite{BM} and standard elliptic regularity results.

\begin{lemma} \label{lem:dani}
Let $k \in C(\D \backslash \{ 0\})$ be a bounded nonpositive continuous function and
let $u \in C^2(\D \backslash \{ 0 \})$ be a real--valued solution to the curvature
  equation (\ref{eq:gauss}) which satisfies (\ref{eq:int}).
Then there exists a real number $\gamma>-1$ and a harmonic function $h$ on
$\D$ such that
$$ u(z)=\gamma \log |z|+h(z) +  v(z) \, , $$
where
$$ v(z)= \frac{1}{2\pi} \iint \limits_{\D}
\log|z-\zeta| \, k(\zeta) \, e^{2
  u(\zeta)} \, d\sigma_{\zeta} \, $$
is the Newton potential of $k(z) \, e^{2 u(z)}$ and $d\sigma_{\zeta}$ denotes
two--dimensional Lebesgue measure w.r.t.~$\zeta$.
\end{lemma}

{\bf Proof.}

(a) \, It turns out that it is more appropriate to work with the 
{\it nonpositive} Green potential
of $k(z) \, e^{2 u(z)}$ instead of its Newton potential $v$.
We therefore let 
$$ g_{\D}(z,\zeta):=-\log \left| \frac{z-\zeta}{1-\overline{\zeta} \, z }
\right| \, , \qquad z, \zeta \in \D\,$$
denote Green's function on $\D$, i.e., $g(z,\zeta) \ge 0$, and let
$$ q(z):=\frac{1}{2 \pi} \iint \limits_{\D} g_{\D}(z,\zeta) \, k(\zeta) \,
e^{2 u(\zeta)} \, d\sigma_{\zeta} \, .$$
Thus $q$ is a well--defined nonpositive continuous function on $\D \backslash \{ 0 \}$ and
  can be decomposed as $q(z)=-v(z)+h_1(z)$, where $v$ is the Newton potential of
$k(z) \, e^{2 u(z)}$ and $h_1$ is a harmonic function in $\D$ (including $z=0$).
Thus $q$ is subharmonic in $\D$ (including $z=0$; we do not know yet whether $q(0)>-\infty$), see \cite[Theorem 3.1.2.]{Ransford}. Now, let
$$ w:=q+u \, . $$
Then $w$ is continuous on $\D \backslash \{ 0\}$ and
$$\Delta w(z)=\Delta q(z)+\Delta u(z)=-k(z) \, e^{2 u(z)}+k(z) \, e^{2 u(z)}\, =0 
\quad \text{ in } \D \backslash\{ 0\}$$
in the sense of distributions, see \cite[p.~74]{Ransford}. By Weyl's lemma,
 $w$ is harmonic in $\D \backslash \{
0 \}$. Therefore
$$ w(z)=\gamma \log |z|+\Re g(z) \, , $$
for some constant $\gamma \in \R$ and some holomorphic function $g : \D
\backslash \{ 0 \} \to \C$.

\medskip

(b) \, We now prove that $g$ is in fact holomorphic on the whole unit disk
$\D$. If not, then $g$ would have a nonremovable singularity at $z=0$, so
$e^{ g}$ would have an essential singularity at $z=0$:
\begin{equation} \label{eq:q}
 e^{g(z)}=\sum \limits_{n=-\infty}^{\infty} b_n z^n \, , \quad b_n\not=0
\text{ for infinitely many n}<0 \, .  
\end{equation}
However,  using Parseval's identity,
\begin{eqnarray*}
  2 \pi \sum \limits_{n=-\infty}^{\infty} |b_n|^2 \int \limits_{0}^{1} r^{2 n+2
  \gamma+1} \, dr &=& \int \limits_{0}^{1} \int \limits_{0}^{2 \pi} 
\left| \sum \limits_{n=-\infty}^{\infty} b_n r^n e^{i n t} \right|^2
\, dt \, r^{2 \gamma+1} \, dr= \iint \limits_{\D} |z|^{2 \gamma} \, \left|
  e^{g(z)} \right|^2 \, dxdy \\ & \le & 
\iint \limits_{\D}   e^{2  \, u(z)}  \, dxdy < \infty \, . 
\end{eqnarray*}
Consequently, $b_n=0$ for all $n< -1-\gamma$, which contradicts (\ref{eq:q}).

\medskip

(c) \, We next show that $\gamma>-1$.

\smallskip

Since $g$ is holomorphic in $\D$, we have
$$u(z)=w(z)-q(z)=\gamma \log |z|+\Re(g(z))-q(z) \ge \gamma \log |z|+\Re
(g(z)) \ge \gamma \log |z|+c$$ 
in $|z|<1/2$ for some real constant $c$, so
$$ e^{2 u(z)} \ge |z|^{2 \gamma} \, |e^{2 g(z)}| \ge |z|^{2 \gamma}  e^{2 c}  \, , \qquad |z|<1/2 \, . 
$$
Thus the assumption (\ref{eq:int}) implies 
$$ \iint \limits_{|z|<1/2} |z|^{2 \gamma} \, dxdy<\infty \, , $$
so $\gamma >-1$.

\medskip

(d) \, Returning to what we have proved in (a), we see that $u(z)=\gamma \log
|z|+h(z)+v(z)$, where $h(z):=\Re (g(z))-h_1(z)$ is harmonic in $\D$ and $v$ is
the Newton potential of $k(z) \, e^{2 u(z)}$.~\hfill{$\blacksquare$}

\begin{lemma}[Br\'ezis--Merle \cite{BM}] \label{lem:brezis}
Let $f \in L^1(\D)$  and let $v(z)$
be the Newton potential of $f$. Then
$$ e^{|v|} \in L^p(\D) \quad \text{ for every } 0<p<+\infty \, . $$
\end{lemma}

{\bf Proof of Theorem \ref{thm:main}.}
By Lemma \ref{lem:dani}, $u(z)=\gamma \log |z|+r(z)$, where $r \in C^2(\D \backslash \{ 0 \})$. Thus, we only need to investigate the regularity properties
of the remainder function $r(z)$ in {\it some} neighborhood of $z=0$, say in $\D_{1/2}=\{z \in \C \, : \, |z|<1/2\}$.
Let  $v$ be the Newton potential of 
$f(z):=k(z) \, e^{2 u(z)}$.  By assumption $f \in L^1(\D) \cap C(\D \backslash \{
0 \})$, so Lemma \ref{lem:brezis} shows  $e^{v} \in L^p(\D)$ for any $0<p<+\infty$.
Using Lemma \ref{lem:dani}  we can write $u(z)=\gamma \log|z|+r(z)$ for some
$\gamma >-1$ with $r=h+v$ where $h$ is harmonic in $\D$. We thus get
$$ f(z)=k(z) \, |z|^{2 \gamma} e^{2 h(z)} e^{2 v(z)} \, . $$
We first consider the case $-1<\gamma<0$. Then $f \in L^q(\D_{1/2})$ for every
 $1<q< -1/\gamma$.  The regularity properties of the Newton potential
show that $v$ belongs to the Sobolev space $W^{2,q}(\D_{1/2})$ for any $1<q <-1/\gamma$, see \cite[Theorem 9.2.1]{Jo}.
Thus the Sobolev embedding theorems yield $r=h+v \in C(\D_{1/2})$ if $-1<\gamma \le
-1/2$ and $r=h+v \in C^1(\D_{1/2})$ if $-1/2<\gamma<0$. The case $\gamma \ge 0$ is
easier, because now $f \in L^q(\D_{1/2})$ for any $1<q<\infty$. Hence $v \in
W^{2,q}(\D_{1/2})$ for any $1<q <\infty$, so $r=h+v \in C^1(\D_{1/2})$. 
If $k$ is locally H\"older
continuous on $\D$ (and $\gamma \ge 0$), then $f$ is locally H\"older
continuous on $\D$, so $r=h+v \in C^2(\D)$.~\hfill{$\blacksquare$}

\bigskip

{\bf Proof of Corollary \ref{cor}.}
We first consider the case $\gamma \ge 0$. Then $u(z)=\log \lambda(z)=\gamma \log
|z|+r(z)$ with $r \in C^1(\D)$ by Theorem \ref{thm:main}, so
$$ z\, \Gamma_{\lambda}(z)= \gamma +2 z \frac{\partial r}{\partial z}(z)\to \gamma
\quad \text{ as } z \to 0 \, .$$
If $k$ is locally H\"older continuous in $\D$, then $r \in C^2(\D)$ and
$$ z^2 \, S_{\lambda}(z)= -\frac{\gamma (2+\gamma)}{2}-2 \gamma z \frac{\partial
  r}{\partial z}(z)+2 z^2\left( \frac{\partial^2 r}{\partial z^2}(z)-\left(
    \frac{\partial r}{\partial z}(z) \right)^2 \right) \overset{z \to
  0}{\longrightarrow} -\frac{\gamma (2+\gamma)}{2}\, . $$
Now let $-1<\gamma < 0$. Choose a positive integer $m$ with $\gamma^*:=\gamma m+m-1 \ge 0$.
Let $\pi_m : \D\backslash \{ 0 \} \to \D \backslash \{ 0 \}$ be the holomorphic
cover projection $\pi_m(z)=z^m$. Then the pullback 
$\lambda^*(z) \, |dz|:=m |z|^{m-1} \lambda(z^m)\, |dz|$ of $\lambda(z) \,
|dz|$ via $\pi_m$ 
induces a conformal metric on $\D \backslash \{ 0\}$ with curvature $- k(z^m)$ 
and $u^*(z):=\log \lambda^*(z)=\gamma^* \log |z|+r^*(z)$.
Moreover, using the transformation rules for the connection and the projective
connection (see \cite[p.~334]{Min97}),
$$ \Gamma_{\lambda^*}(z)=m \Gamma_{\lambda}(z^m)z^{m-1}+\frac{m-1}{z} \, ,
\qquad S_{\lambda^*}(z)=S_{\lambda}(z^m) m^2 z^{2(m-1)}-\frac{m^2-1}{2\,
  z^2}\, ,
$$
we get
$$ \lim \limits_{z \to 0} z \Gamma_{\lambda}(z)
=\lim \limits_{z \to 0} z^m \Gamma_{\lambda}(z^m)=
\frac{\gamma^*-(m-1)}{m}=\gamma \,  $$
and, if $k$ is locally H\"older continuous on $\D$, then
$$ \lim \limits_{z \to 0} z^2 S_{\lambda}(z)=\lim \limits_{z \to 0} z^{2m}
S_{\lambda}(z^m)=\frac{-\gamma^* (2+\gamma^*)+m^2-1}{2\, m^2}=-\frac{\gamma
  (2+\gamma)}{2} \, . 
$$
\hfill{$\blacksquare$}

\end{document}